\documentclass[english]{amsart}

\usepackage{babel}
\usepackage{amstext}
\usepackage{amsmath}
\usepackage{amsfonts}
\usepackage{latexsym}
\usepackage{ifthen}
\usepackage{xypic}
\xyoption{all}
\newcommand\sI{{\mathcal I}}

\newcommand\sL{{\mathcal L}}
\newcommand\sO{{\mathcal O}}

\newcommand\sX{{\mathcal X}}
\newcommand\sC{{\mathcal C}}

\newcommand\bR{{\mathbb R}}
\newcommand\bZ{{\mathbb Z}}
\newcommand\bC{{\mathbb C}}
\newcommand\bQ{{\mathbb Q}}
\newcommand\bN{{\mathbb N}}

\newcommand\bP{{\mathbb P}}
\newcommand\rh{{\dasharrow}}

\newcounter{lemma}

\newtheorem{lemma1}[lemma]{\setcounter{equation}{0}}

\newenvironment{lemma}{\begin{lemma1}{\bf Lemma.}}{\end{lemma1}}

\newenvironment{theorem}{\begin{lemma1}{\bf Theorem.}}{\end{lemma1}}

\newenvironment{proposition}{\begin{lemma1}{\bf Proposition.}}{\end{lemma1}}

\newenvironment{corollary}{\begin{lemma1}{\bf Corollary.}}{\end{lemma1}}

\newenvironment{definition}{\begin{lemma1}{\bf Definition.}}{\end{lemma1}}

\newenvironment{setup}{\begin{lemma1}{\bf Setup.}}{\end{lemma1}}

\newenvironment{conjecture}{\begin{lemma1}{\bf Conjecture.}}{\end{lemma1}}

\newenvironment{Induction Step}{\begin{lemma1}{\bf Induction Step.}}{\end{lemma1}}
\newenvironment{Proof of Theorem 1.2}{\begin{lemma1}{\bf Proof of Theorem 1.2.}}{\end{lemma1}}
\newenvironment{can}{\begin{lemma1}{\bf Conventions and Notations}}{\end{lemma1}}

\title {Non-algebraic hyperk\"ahler manifolds} \author{Fr\'ed\'eric Campana, Keiji Oguiso and Thomas Peternell}
\date{\today}

\begin{document}

\maketitle

\tableofcontents

%\vspace*{-0.5in}
\section{Introduction}

Let $X$ be a compact simply connected complex irreducible symplectic K\"ahler manifold of dimension $2n$ (a hyperk\"ahler manifold, for
short), that is a
compact simply connected K\"ahler manifold
admitting a holomorphic $2-$form $\sigma_X$ which is of maximal rank at every point such that $H^0(X, \Omega_X^2) = \bC \sigma_X$ 
(hence $\wedge^n \sigma_X$ is a $2n-$form without zeroes). By Fujiki \cite{Fu83}, based on Bogomolov's unobstructedness theorem
\cite{Bo78}, both projective and non-projective hyperk\"ahler 
manifolds are dense in the Kuranishi space of $X$.

We are concerned with non-algebraic hyperk\"ahler manifolds, particularly with their algebraic dimensions $a(X) <  2n$ and their algebraic
reductions $ f: X \dasharrow B$, which are unique up to bimeromorphic modification of $B$.

A central role in the theory is played by the Beauville form 
$$q_X: H^2(X,\bZ) \times H^2(X,\bZ) \to \bZ $$
of signature $(3, 0, b_2(X)-3)$ \cite{Be83}. It induces a symmetric bilinear form on the N\'eron-Severi group $NS(X)$,
and the signature is either $(1, 0, \rho(X)-1)$ in which case we say that $NS(X)$ is hyperbolic,
$(0,1, \rho(X)-1)$ $(NS(X)$ is parabolic), or $(0,0, \rho(X))$ ($NS(X)$ is elliptic). Now $X$ is projective (equivalently Moishezon
\cite{Mo66}) if and only if $NS(X)$ is hyperbolic (Huybrechts \cite{Hu99}), so we shall be interested in the parabolic case and the elliptic
case. \\
The Beauville form can be seen as a natural higher dimensional version of the intersection form on a compact complex surface $S$. 
Here we have
$a(S) = 0$, $1$, $2$ according to $NS(S)$ being elliptic, parabolic, hyperbolic. In addition, if
$a(S) = 1$, then we have a holomorphic algebraic reduction $f : S \longrightarrow C$ whose general fiber is an elliptic curve 
\cite{BHPV04}. Of great importance in the theory of hyperk\"ahler manifolds is the following
 example due to Beauville \cite{Be83} and Fujiki \cite{Fu83}.
Let $S$ be a K3 surface. Then $S^{[n]}$, the Hilbert scheme of $n$ points on $S$, is a hyperk\"ahler manifold of dimension $2n$. 
We have $a(S^{[n]}) = 0$, $n$, $2n$ according to $a(S) = 0$, $1$, $2$. In addition, when $a(S) = 1$, the algebraic reduction map 
$S \to \bP^1$ induces a natural morphism $S^{[n]} \longrightarrow \bP^n$. This is the algebraic reduction of $S^{[n]}$ and
it is also Lagrangian. \\
This motivates the following \cite{Og07}

\begin{conjecture}\label{conjecture:hkc}
Let $X$ be a hyperk\"ahler manifold of dimension $2n.$ Then its algebraic dimension takes only the values $0,n,2n.$
Moreover, if $a(X) = n$, then the algebraic reduction has a holomorphic model $ f: X \longrightarrow B$ with $B$ a normal projective
variety of dimension $n$. Finally
$f$ is Lagrangian, that is $\sigma_X \vert F \equiv 0$ for a general fiber of $f$. 
\end{conjecture}

The last statement actually follows from the previous by \cite{Ma99}.

The aim of this paper is to establish the following (partial) answers to Conjecture 1.1. 

\begin{theorem}\label{4dim} Conjecture 1.1 holds in dimension 4. 
\end{theorem} 

In higher dimensions we can determine the algebraic dimension up to the existence of minimal models of K\"ahler spaces with algebraic dimenion
and Kodaira dimension $0:$ 

\begin{theorem}\label{theorem:main0} Let $X$ be a hyperk\"ahler manifold of dimension $2n.$ Then Conjecture 1.1 holds
provided any compact K\"ahler manifold $Y$ with $\dim Y \leq 2n-1,$ $a(Y) = \kappa (Y) = 0$ has a minimal model. 
\end{theorem}

A general conjecture from minimal model theory says that every compact K\"ahler manifold of non-negative Kodaira dimension 
should have a minimal model. \\
Recall that a meromorphic map between compact
varieties is {\it almost holomorphic}
if the general fiber is compact, i.e., does not meet the set of indeterminacies. Then we can state more precisely:  

\begin{theorem}\label{theorem:main} Let  $X$ be a hyperk\"ahler manifold of dimension $2n.$ Then:
\begin{enumerate}
\item If $NS(X)$ is elliptic, then $a(X) = 0$.
\item If $NS(X)$ is parabolic, then $0 \le a(X) \le n = \dim\, X/2$.
\item Assume that $NS(X)$ is parabolic and $a(X) > 0.$ Then one can choose an algebraic reduction of one of the following two forms:

    (i) $f : X \longrightarrow B$ is holomorphic Lagrangian, in particular, $a(X) = n$, or

    (ii) $f : X \dasharrow B$ is not almost holomorphic
    and the general fiber $X_b$ ($b \in B$) is isotypically semi-simple, in particular, $a(X_b) = 0$ (see section 2
for the definition of the term ``isotypically semisimple'').
\item Assume that any compact isotypically semi-simple K\"ahler manifold $Y$ of $\dim\, Y \le 2n-1$, of algebraic dimension $a(Y) = 0$ and of Kodaira dimension 
$\kappa(Y) = 0$ and with effective canonical divisor $K_Y$,  has a minimal model. Then Conjecture 1.1 holds.
\end{enumerate}
\end{theorem}

\vskip .2cm \noindent One might speculate that a hyperk\"ahler manifold $X$ of dimension $2n$ has algebraic dimension
$a(X) = n$ if and only if $NS(X)$ is parabolic. This in turn would be a consequence of a potential semi-ampleness of any nef
line bundle on a hyperk\"ahler manifold. 
\vskip .2cm \noindent Parts (1) and (2) of Theorem 1.4 will be proved in section 3; parts (3) and (4) in section 4 and 
Theorem 1.2 finally in section 5. All these sections make essential use of section 2, which contains structure results on
meromorphic fibrations on 
compact K\"ahler manifolds, in particular on those manifolds admitting a unique holomorphic 2-form which additionally is
generically non-degenerate. The final section gives some results on nef line bundles on hyperk\"ahler manifolds.

\vskip .5cm \noindent
{\it Acknowledgement.} Our collaboration has been made possible by
the priority program {\it  ``Global methods in complex
geometry''} and the research group {\it ``Algebraic surfaces and compact complex manifolds} of the Deutsche Forschungsgemeinschaft, which we
gratefully acknowledge. In the final stage we profited very much from a joint stay at the Korea Institute of Advanced Study; we would
like to thank Jun-Muk Hwang for the invitation and the excellent working conditions.

%\vfill \eject

\tableofcontents

\section{Fibrations on generalized Hyperk\"ahler Manifolds}\label{min}
\setcounter{lemma}{0}
In this section we prove some general structure theorems on generalised hyperk\"ahler manifolds.

\begin{can} {\rm
(1) By $X,X',X^{''}, \ldots$, we denote $n$-dimensional compact irreducible complex spaces which are bimeromorphic to compact 
K\" ahler manifolds.
The algebraic dimension ([Ue75]) is denoted by $a(X)$. We have the algebraic reduction map (only defined up to obvious  
bimeromorphic equivalence) 
$$f:X \dasharrow B.$$ We always take $B$ to be normal projective and often we will choose $B$ smooth. 
\vskip .2cm \noindent
(2) A {\it fibration} $f:X \dasharrow Y$ is a dominant meromorphic map with connected fibres.   
The fibration $f$  is said to be {\it almost holomorphic}  if its generic fibre does not meet the indeterminacy locus of $f$.  
If $Y$is not uniruled, 
then any fibration $f: X \dasharrow Y$ is automatically almost holomorphic. \\ 
The fibration $f$ is said to be {\it trivial} if $\dim Y = 0$ or $\dim Y = n$. 
\vskip .2cm \noindent 
(3) A point $y$ in $Y$ is said to be {\it general} if it lies outside of a countable union of (suitable) proper closed analytic subsets of $Y$. 
We denote by $X_y$ the (Chow-theoretic) fibre of $f$ over a generic $y\in Y$. 
\vskip .2cm \noindent 
(4) Recall ([Fu83]) that a compact K\"ahler manifold $X$ is said to be {\it simple} if $X$ is not covered by positive-dimensional 
irreducible compact proper analytic subsets.
If $X$ is simple, then necessarily $a(X)=0$. Moreover either $q(X)=0$ or there $X$ is (bimeromorphically) covered by a simple torus, i.e., 
there exists a simple non-algebraic 
complex torus $T$ and a meromorphic dominant map $u:T \rh X$. \\ 
Remark that if $X$ is a complex torus with $a(X)=0$, $X$ is simple in the sense above if and only if it is simple in the classical sense, 
i.e., $X$ has no nontrivial complex subtorus. However, if $X$ is an abelian variety, it is not simple in the above sense, 
but can be simple in the 
classical sense. 
%By an abuse of notation, we shall say that an {\it abelian variety} is simple (or A-simple, to avoid confusion) if it is so in the classical sense. 
\vskip .2cm \noindent 
(5) Two complex spaces $X$ and $X'$ are {\it commensurable} if there exist a complex space $X^{''}$ and generically finite surjective
holomorphic maps $X^{''}\to X$ and 
$X^{''} \to X'$. This is (easily seen to be) an equivalence relation. \\ 
Notice that two {\it projective} varieties are commensurable if and only 
if they have the same 
dimension. But for non algebraic $X$, this equivalence relation is very restrictive. If $X$ is simple, and if $X'$ is commensurable to $X$, then $X'$ is 
simple, too. 
Thus either $q(X')=0$ for any $X'$ commensurable to $X$, or $X$ is covered by a simple torus. 
\vskip .2cm \noindent 
(6) We say (after \cite{Fu83}) that $X$ is {\it semi-simple} if it is commensurable to a product of simple manifolds, and that $X$ is {\it isotypically semi-simple } 
if it is commensurable to a product $S^k$ for some simple $S$ and some $k>0$. 
Remark that if $f:X\dasharrow S$ is a fibration with $S$ semi-simple, then $f$ is almost holomorphic, since $S$ is not uniruled.\\ 
%In completely the same way we define the notion of ``isotypical A-semi-simplicity''.
} 
\end{can}

\begin{definition} Let $f:X \dasharrow  B$ be a fibration from a compact (connected) K\"ahler manifold $X$.\\
We say that  $f=h\circ g$ is a factorisation of
$f$ if $g:X\rh S$ and $h:S\rh B$ are fibrations with $f=h\circ g$.\\
This factorisation is said to be trivial if $\dim S = \dim X$ or $\dim S = \dim B.$ 
The fibration $f$ is said to be {\it minimal} if
any factorisation of $f$ is trivial.
The variety itself $X$ is minimal if the constant fibration $X \to \{pt\}$ is minimal.
\end{definition}

The following easy observation is essential:

\begin{theorem} Let $X$ be a compact K\"ahler manifold of dimension $2n$ and suppose that $h^{2,0}(X) = 1$, and the corresponding
holomorphic $2-$form $\sigma$ (which is unique up to a scalar) satisfies $\sigma^n\neq 0$.
Then the following assertions hold.
\begin{enumerate}
\item If $f: X \rh Y$ is a fibration with $\dim Y < \dim X$, then $Y$ is Moishezon.
\item The algebraic reduction $f: X \rh B$ is minimal.
\item If $a(X) = 0$, then $X$ is minimal.
\item In particular hyperk\"ahler manifolds $X$ with $a(X) = 0$ (see section 3 below) are minimal. 
\end{enumerate}
\end{theorem}

\proof Only (1) needs to be proved; the other statements are trivial consequences.
Suppose that $f:X \rh  Y$ is a nontrivial fibration with $Y$ a non-algebraic manifold (we may assume $Y$ smooth and K\"ahler). 
Then $h^{2,0}(Y)>0$, by \cite{Ko54}.
Any non-zero holomorphic $2-$form on $Y$ lifts to the holomorphic $2-$form $\sigma$  on $X$, which is generically of maximal rank, 
so that $\dim Y = \dim X$. 
\qed

\vskip .2cm \noindent
The main result of this section is:

\begin{theorem}\label{minfib} Let $X$ be a compact K\"ahler manifold, $f:X\rh B$ be a minimal fibration. Suppose  $\dim X > \dim B.$
Let $X_b$ be a general fibre of $f$. Then
\begin{enumerate}
\item $X_b$ is either Moishezon or isotypically semi-simple, in which case $a(X_b) = 0.$
\item Furthermore, if $X$ is not projective and if $B$ and $X_b$ are Moishezon, then $f$ is almost holomorphic and $X_b$ is an
abelian variety.
\end{enumerate}
\end{theorem}

\ref{minfib} will be proved at the end of this section.

\vskip .2cm \noindent
Applying 2.3(2), we get:

\begin{corollary}\label{algred} Let $X$ be a compact K\"ahler manifold of dimension $2n$ with $h^{2,0}(X)=1$,
carrying a holomorphic two-form  $\sigma$ such that $\sigma^n\neq 0$.
Assume $X$ is nonprojective. Then the following assertions hold.
\begin{enumerate}
\item Let $f:X \dasharrow B$ be the algebraic reduction and $X_b$ be a general fibre of $f$. Then either
\subitem 2a. $X_b$ is isotypically semi-simple or
{\subitem 2b. $f$ is almost holomorphic, and $X_b$ is an abelian variety.}
\item In particular if  $a(X)=0$, then $X$ is isotypically semi-simple.
\end{enumerate}
\end{corollary}

\begin{corollary} Let $X$ be a simply-connected compact K\"ahler manifold of dimension $2n$ with $h^{2,0}(X)=1$,
carrying a holomorphic two-form  $\sigma$ such that $\sigma^n\neq 0$. Assume moreover that $X$ does not contain any effective divisor. 
Then $X$ is simple.

In particular, any simply connected irreducible hyperk\" ahler manifold without effective divisors is simple 
and so does the general member of the Kuransihi family of any 
hyperk\" ahler manifold without effective divisors is simple.
\end{corollary}

\proof Since $X$ has no effective divisors, we have $a(X) = 0.$ Thus $X$ is isotypically semi-simple by \ref{algred}.
Specifically there exist generically finite meromorphic maps $u:Z\to X$ and $v:Z\to S^k$,
with $Z$ smooth, and $S$ simple. Our claim comes down to prove that $u$ is bimeromorphic and that $k=1$.
Since $X$ has no divisor, $u$ is unramified, hence bimeromorphic, $X$ being simply-connected.
Thus $h^{2,0}(X)=h^{2,0}(Z)=1$. Since $S$ is non-algebraic, one has $h^{2,0}(S)\geq 1$, hence necessarily $k = 1.$
\qed

\vskip .2cm \noindent
The assumption that $X$ does not contain any divisor cannot be removed in (2.6). In fact, the hyperk\"ahler 4-fold $S^{[2]}$
with $S$ a K3 surface, $a(S) = 0$ is not simple but $a(S^{[2]}) = 0.$ 

The rest of the section is devoted to the proof of \ref{minfib}.
We shall need the following two elementary lemmas, which are relative versions of results similar
to those in [Fu83] in a simplified form. We first recall some notions needed in the proof.

A {\it covering family} of $X$ will be a compact irreducible analytic subset $S\subset \sC(X)$ of the Chow
(or cycle, or Barlet)-space $\sC(X)$
(see [Ba75]) of $X$, such that if $Z\subset S\times X$ is its incidence graph, with natural projections $p:Z\to X$ and $q:Z\to S$,
then $p$ is surjective,
and the generic fibre of $q$ is irreducible. In other words,
$X$ is covered by the (generically irreducible) cycles $Z_s$, $s\in S $. We call $m = \dim Z_s$ the dimension of the family $S$. \\
If $f: X \dasharrow Y$ is a fibration, we denote by $\sC(X/Y)$ the closed analytic subset of $\sC(X)$ consisting of those points $s\in \sC(X)$ such that
the corresponding analytic compact pure-dimensional cycle $Z_s$ of $X$ has support contained in one fibre $X_y$ of $f$.
If $S \subset \sC(X/Y)$ is a covering family of $X$, the map $f_*$ sending $s$ to $y=f(Z_s)$ is a meromorphic dominant map $f_*: S \dasharrow Y$.
\vskip .2cm \noindent

\

\begin{lemma}\label{covf} Let $X$ be a compact K\"ahler manifold and $f: X \rh B$ be any fibration with $a(X_b) = 0.$
Let $S\subset \sC(X/B)$ be a nontrivial
covering family of $X$ over $B$.
Assume that $\dim Z_s = m$ is
maximal among the dimensions of nontrivial covering families of $X$ over $B$. Then
\begin{enumerate}
\item $\dim Z = \dim X$ (hence only finitely many of the $Z_s's$ pass through the generic point of $X$).
\item $S_b$ is simple, and no proper closed analytic subset of $S_b$ is a covering family of $X_b$.
\item $S_b$ is the union of finitely many irreducible components of $\sC(X_b)$.
\end{enumerate}
\end{lemma}

\proof (1) Assume that $\dim Z > \dim X$. Then also $\dim Z_b > \dim X_b$ so that we may assume $\dim B = 0.$
The fibres of $p:Z\to X$ are Moishezon by [Ca81]. Thus we can find a covering family of $S$
by curves $(C_v)_{v\in V}$.
For general $v \in V$ we define
$$W_v:=p(q^{-1}(C_v)).$$
This is an irreducible compact analytic subset of $X$ and defines a covering
family $(W_v)_{v\in V}$ of $X$ with $\dim W_v = m+1$. In order to show that this family is non-trivial, we prove that $m+1 <  \dim X.$
In fact, if $\dim X \leq m+1,$ then $\dim X = m+1$ so that the $Z_s$ are divisors in $X.$ This contradicts $a(X) = 0.$\\
Thus the family $(W_v)$ is non-trivial, contradicting the maximality of $m.$ \\
(2) The same argument shows that $S_b$ is simple. In fact, if $S_b$ were not simple, then in every fiber we find a covering
family of proper subvarieties and in total
we obtain a nontrivial covering family $(C_v)_{v\in V}$ of $S$ and define $W_v$ as above.
By the maximality of $m$, we must have $W_v=X$ for all $v\in V$. But then $(Z_s)_{s\in C_v}$ is a nontrivial covering family of $X$.
Let $Z':=q^{-1}(C_v)$ be the graph of this covering family. Then $\dim Z' = \dim X = \dim Z$. Thus by irreducibility we obtain $Z'=Z$ and
$C_v=S$, a contradiction. \\
(3) The third assertion is an obvious consequence of the second.
\qed

\begin{lemma}\label{subsss} Let $k \in \bN$ and let $S_j$ be simple manifolds for $1 \leq j \leq k$. Put
$$ S = S_1\times \dots \times S_k$$
with projections $p_j: S \to S_j.$ More generally, for a given subset  $J=\lbrace j_1,\dots ,j_h\rbrace \subset \lbrace 1,2,\dots, k\rbrace$,
let
$$p_J: S \to S_J = S_{j_1} \times \ldots \times S_{j_h} $$
be the projection.
Let $Y \subset S$ be an irreducible
compact analytic subset such that $p_j(Y) = S_j$ for all $j.$ \\
There exists $J$ such that $p_J:Y\to S_J$ is surjective and generically finite. In particular, $Y$ is commensurable to $S_J$ and is therefore semi-simple.
In particular, if $S_j \simeq S_k$ for all $j,k$, then $Y$ is isotypically semi-simple.
\end{lemma}

\proof Let $K = \lbrace 1,\dots, k-1\rbrace$. If $p_K(Y) \ne S_K,$ we proceed by induction on $k.$
Thus we may assume that $p_K(Y) = S_K$ for $K=\lbrace 1,\dots, k-1\rbrace$.
If $p_K:Y\to S_K$ is not generically finite, let $S'_K$ be its Stein factorisation with map $p'_K: Y \to S'_K$, and define a meromorphic map
$$\varphi:S'_K\to \sC(S_k)$$
by sending a general $s\in S'_K$ to $p_k({p'_K}^{-1}(s))$. The image of $\varphi$ gives a covering family of $S_k$.
Because $S_k$ is simple, we must have $\varphi(s)=S_k$ for all $s$. Thus $Y=S$ (in which case we take $J = \lbrace 1, \ldots, k \rbrace$).
\qed

\vskip .2cm \noindent 
\proof {\it of \ref{minfib}} Let $a_f:X\rh Y, \ h:Y\rh B$ with $f=h\circ a_f$ be the relative algebraic reduction \cite{Ca81},
\cite{Fu83} of $f$
(so that the restriction
of $a_f$ to $X_b$ is the algebraic reduction of $X_b$, for $b$ general in $B$). Since $f$ is minimal, either $Y=X$ up to 
bimeromorphic equivalence and $X_b$ is Moishezon, or $Y=B$ up to bimeromorphic equivalence
and $a(X_b)=0$.
\vskip .2cm \noindent
(1) In the first step we assume that $a(X_b)=0$ and need to show that $X_b$ is isotypically semi-simple. If $X_b$ is simple, we are already 
done. 
If $X_b$ is not simple, let $S\subset \sC(X/B)$ be a nontrivial covering family of $X$ with
$m = \dim Z_s$ maximal. By (\ref{covf}), $S_b$ is simple and $p:Z\to X$ is generically finite onto $X$. Let $\delta$ be the degree of $p$ and
$$\varphi:X\rh Sym^{\delta}(S/B)$$
be the meromorphic map sending a general $x\in X$ to $q_*(p^{-1}(x))$ (here $Sym^{\delta}(S/B)$ denotes the subspace of
$Sym^{\delta}(S)$ consisting of $\delta$-tuples of $S$ contained in some fibre of $S$ over $B$. We adopt a similar convention for $(S^{\delta}/B)$). \\
Since $f$ is minimal, this map is generically finite onto its image $Y_0\subset Sym^{\delta}(S/B)$. Let $Y\subset (S^{\delta}/B)$ be a main
component of the inverse image of $Y_0$ under the natural map $(S^{\delta}/B)\to Sym^{\delta}(S/B)$. Then $Y$ maps surjectively onto $S$ under all
projections from $(S^{\delta}/B)$ to $S$ (otherwise there would exist some irreducible proper compact analytic subset $S'\subset S$ parametrising
a covering family of $X$, contradicting (\ref{covf})). From (\ref{subsss}) we conclude that $Y$, and hence so $X$,  is commensurable to $(S^k/B)$ for
some $k\leq {\delta}$. The first assertion of (\ref{minfib}) is thus established.
\vskip .2cm 
(2) Thus we now assume that $X_b$ is Moishezon. By the minimality assumption, it follows that $f$ is the algebraic reduction of $X$.
From \cite {Ca81} we deduce that $f$ is almost holomorphic. Moreover by \cite{Fu83} the general fiber $X_b$ is abelian or a unirational
manifold. \\
If $\kappa (X) \geq 0 $ - and this is sufficient for all our applications - $X_b$ cannot be uniruled and we conclude. But if $X_b$ is unirational,
then \cite{Fu83}, Prop.2.5 implies that $X$ is projective, which is excluded
by assumption. 
\qed
\vskip .2cm \noindent 
We now consider the restriction of holomorphic $2$-forms to fibers.

\begin{corollary} Let $X$ be a compact K\"ahler manifold and $f:X\rh B$ be a fibration. Assume
that the restriction of any holomorphic $2$-form on $X$ to the fiber $X_b$
for generic $b\in B$ vanishes. Then
\begin{enumerate}
\item $X_b$ is Moishezon for all  $b$.
\item If $f$ is the algebraic reduction of $X$ and a minimal fibration, then $f$ is almost holomorphic, and $X_b$ is an abelian variety.
\item In particular, suppose that $h^{2,0}(X)= 1$, given by a
holomorphic $2-$form $\sigma$ such that $\sigma^n\neq 0$. Then the algebraic reduction $f$
is almost holomorphic Lagrangian provided $\sigma$ restricted to $X_b$ vanishes.
\end{enumerate}
\end{corollary}

\proof  The first statement is a lemma due to C. Voisin (see [Ca04]). The second follows from \ref{minfib}. 
\qed

\section{Basics on Hyperk\"ahler Manifolds and first results}
\setcounter{lemma}{0}

We begin by fixing some notations.  For the rest of the paper we fix an irreducible simply connected compact complex symplectic 
K\"ahler manifold $X$ of dimension $2n,$ that is a simply connected compact K\"ahler manifold
admitting a holomorphic $2-$form $\sigma$ which is of maximal rank at every point (hence $\sigma^{2n}$ is a $n-$form without zeroes), such that 
$H^0(X,\Omega_X^2) = \bC \sigma. $ 
We say for short
that $X$ is a {\it hyperk\"ahler manifold}. \\
The non-degenerate symmetric bilinear form, constructed by Beauville \cite{Be83}, will be denoted
$$q = q_X: H^2(X,\bZ) \times H^2(X,\bZ) \to \bZ.$$
We shall use the shorthand $q(a) = q(a,a).$
The Picard number of $X$ is denoted by $\rho.$ Then the signature of $q_X$ on the N\'eron-Severi group
$$NS(X) = H^{1,1}(X) \cap H^2(X,\bZ)$$
is one of the following.
\begin{itemize}
\item $(1,0,\rho-1)$ (hyperbolic case);
\item $(0,1,\rho-1)$ (parabolic case);
\item $(0,0,\rho)$ (elliptic case).
\end{itemize} 
Moreover $X$ is projective if and only if $NS(X)$ is hyperbolic \cite{Hu99}, so we are interested in the parabolic case and the elliptic case.

\begin{theorem}\label{theorem:sign}
\begin{enumerate}
\item If $NS(X)$ is elliptic, then $a(X) = 0$.
\item If $0 < a(X) < 2n$, then $NS(X)$ is parabolic.
Let  $\ell \in NS(X)$ be the unique primitive isotropic vector of $NS(X)$ with $q_{X}(\ell, \eta) > 0$ for a K\"ahler class $\eta$. Then there is a line bundle
$L$ whose linear system defines the algebraic reduction, such that $c_1(L) \in \bZ_{>0} \ell.$
In particular $q_X(L) = 0.$
\end{enumerate}
\end{theorem}
\proof Let $f: X \dasharrow B$ be the algebraic reduction with $B$ normal projective.
If $\dim\, B > 0$, i.e., $a(X) > 0,$ then there is a line bundle $L$ with $D_1, D_2 \in \vert L \vert$ such that 
$D_1$ and $D_2$ have no common component. By definition of the
Beauville form, one has (up to positive constant multiple):
$$q_X(L, L) = \int_X c_1(D_1)c_1(D_2)(\sigma \wedge \overline{\sigma})^{n-1} = \int_{D_1 \cap D_2} (\sigma \wedge 
\overline{\sigma})^{n-1} \ge 0\, .$$
Thus $NS(X)$ is not elliptic if $a(X) > 0$. This proves (1). Since
$X$ is projective if $NS(X)$ is hyperbolic, we have from the same inequality that $q_X(L) = 0$. Hence $NS(X)$ is parabolic
if $0 < a(X) < 2n$. Moreover, $q_X(L, \eta) > 0$ by the shape of Beauville form above. This proves (2).
\qed

\begin{setup} {\rm (1)
We shall assume that $0 < a(X) < 2n.$ Thus $X$ is not projective and $NS(X)$ is parabolic.  
We consider ``the'' algebraic reduction
$$ f: X \dasharrow B. $$
From the previous section we recall that the general fiber is isotypically semi-simple or that 
$f$ is almost holomorphic and the general fiber is abelian. \\ 
We always take $B$ to be normal {\it projective} and often we will choose $B$ smooth, too.
Let
$$ \pi: \tilde X \to X $$
be a resolution of indeterminacies of $f$ so that the induced map $\tilde f: \tilde X \to B$ is holomorphic.
\vskip .2cm \noindent
(2) We fix a very ample line bundle $A$ on $B$ and set
$$ L = \pi_*(\tilde f^*(A))^{**}. $$
This is a holomorphic line bundle on $X$ and we find an effective divisor $E$ on $\tilde X$ such that
$$ \pi^*(L) = \tilde f^*(A) + E.$$
We set $\tilde L = \pi^*(L).$
\vskip .2cm \noindent
(3) In all what follows $\eta $ will always denote a K\"ahler form on $X.$ We set $\tilde \eta = \pi^*(\eta).$
\vskip .2cm \noindent}
\end{setup}

By the results of section 2 we may already state:

\begin{proposition} If the algebraic reduction $f: X \to B$ is holomorphic (with $B$ projective and $\dim B > 0$)
then $f$ is Lagrangian, in particular all smooth fibers are abelian. 
\end{proposition}

\proof By \cite{Ma99}, $f$ is Lagrangian and $\dim B = n$ - his argument works in the K\"ahler case as well. Then there is no holomorphic $2-$form with
non-zero restriction to the general fiber. Therefore we conclude by (2.9). 
\qed

\begin{theorem}\label{theorem:kaehler}
Assume that $0 < a(X) < 2n$. Then $c_1(L) \in \overline{\mathcal K(X)}$, the closure of the K\"ahler cone, i.e., $L$ is (analytically)
nef.
Moreover, $(L.C) = 0$ for all curves
$C \subset X$.
\end{theorem}

\proof Let $\overline{\mathcal P(X)} \subset H^{1,1}(X, \bR)$ be the closure of the positive cone of $X$.
By Theorem \ref{theorem:sign} (2),
$c_1(L) \in  \overline{\mathcal P(X)}$. Thus by \cite{Hu99}
$$c_1(L) \in \overline{\mathcal K(X)}$$ if $L \cdot C \ge 0$ for all curves
$C \subset X$. As, the Beauville form $q_X$  is non-degenerate and defined over $H^2(X, \bQ)$, the map
$x \mapsto q_X(*, x)$ gives an isomorphism
$$\iota : H^2(X, \bQ) \simeq H^{4n-2}(X, \bQ).$$
Here we identify $H^{4n-2}(X, \bQ)$ with $H^2(X, \bQ)^{*}$ by the intersection pairing. Moreover, by the shape of the
Beauville form, $\iota$ induces an isomorphism
$$\iota : H^{1,1}(X, \bQ) \simeq H^{2n-1,2n-1}(X, \bQ).$$ Since $[C] \in H^{2n-1,2n-1}(X, \bQ)$, there is an element
$\alpha \in H^{1,1}(X, \bQ)$ such that $\iota(\alpha) = [C]$. Hence
$$L \cdot C = q_X(L, \alpha) = 0$$
by Theorem \ref{theorem:sign} (2). This proves our claim.
\qed

\begin{lemma}
\begin{enumerate}
\item $L^n \cdot \eta^n > 0 $.
\item $L^{n+1} \alpha_1 \cdot \ldots  \alpha_{n-1} = 0 $ for all $\alpha_i \in H^{1,1}(X).$
\item $L^{n+1} = 0$ in $H^{2n+2}(X,\bR),$ hence the numerical dimension $\nu(L) = n.$ 
\item For all $a,b \geq 0$ with $a+b > n$ we have
$$ \tilde f^*(A)^a \cdot \tilde L^b \cdot \tilde \eta^{2n-(a+b)} = 0.$$
\item For all $k \geq 0$ we have
$$ \tilde f^*(A)^k \cdot \tilde L^{2n-k} = 0.$$
\end{enumerate}
\end{lemma}

\proof Let $\alpha \in H^{1,1}(X).$ Then we have via $q(L) = 0$, using Fujiki's relation \cite{Fu87}
$$(L + t \alpha)^{2n} = c q(L + t \alpha)^n = c (2tq(L,\alpha) + t^2q(\alpha))^n$$
with $c > 0.$ Comparing coefficients, this shows
$$ L^{2n-k} \cdot \alpha^k = 0$$
for $n > k \geq 1,$ hence (2) by polarization. \\
The first claim is just (up to a positive multiple) $ L^n \cdot \eta^n  = c q(L,\eta)^{n} > 0$ by Theorem 3.1 (2).\\
For the proof of (3) notice that  $c_1(L)^{n+1} $ is represented by a closed positive current $T$ of bidimension $(n+1,n+1)$ 
(approximate e.g. $c_1(L)$ by K\"ahler forms).
Choosing $\alpha_i $ in (2) to be represented by K\"ahler forms $\omega_i$, we conclude that
$$ 0 = L^{n+1} \cdot \alpha_1  \ldots \cdot \alpha_{n-1} = T(\omega_1 \wedge \ldots \wedge \omega_{n-1}).$$
Thus $T = 0$ and therefore $c_1(L)^{n+1} = 0$ in $H^2(X,\bR).$ Hence $\nu(L) = n.$ \\
Let $c > n$ and $a+b = c.$ By (2) we have
$$\tilde L^c \cdot \tilde \eta^{2n-c} = 0.$$
Hence $$ (\tilde f^*(A) + E) \cdot \tilde L^{c-1} \cdot \tilde \eta^{2n-c} = 0.$$
Since $\tilde L$ is nef, this gives
$$ \tilde f^*(A) \cdot \tilde L^{c-1} \cdot \tilde \eta^{2n-c} = E \cdot  \tilde L^{c-1} \cdot \tilde \eta^{2n-c} = 0.$$
Continuing in this way, we obtain
$$ \tilde f^*(A)^a \cdot \tilde L^{c-a} \cdot \tilde \eta^{2n-c} = 0$$
proving (4). \\
Claim (5) is the following special case of (4): $a = k$ and $b = 2n-k.$
\qed

\begin{theorem}\label{corollary:geomineq}
Let $X$ be a non-algebraic hyperk\"ahler manifold of dimension $2n. $ Then
$a(X) \le n$.
\end{theorem}
\proof
Recall that $a(X) = \dim\, B$ and suppose $k = \dim\, B > n$. Then we can take $a = k$ and $b = 0$ in Lemma 3.5, i.e., 
$(\tilde{f}^*A)^k\tilde{\eta}^{2n-k} = 0$. The class $(\tilde f^*A)^k$ is represented by a positive 
multiple of general fiber $\tilde{F}$ of $\tilde{f}$.
So, if we put $F = \pi_{*}\tilde{F}$, then $F$ is a $2n-k$-dimensional non-zero effective cycle on $X$, but
$$F\eta^{2n-k} = \tilde{F}\tilde{\eta}^{2n-k} = 0\, ,$$
a contradiction.
\qed

\begin{theorem} Let $X$ be a hyperk\"ahler manifold of dimension $2n.$ Suppose $a(X) = n$. Then any nef line bundle 
$D$ on $X$ is semi-ample. In particular its algebraic reduction can be taken holomorphic.
\end{theorem}

\proof Since  $NS(X) = \bZ \ell \oplus V$ where $q_X$ is negative definite on $V,$ the line bundle 
$D = L$ up to a multiple. By (3.5) we know $\nu(L) = n,$ on the other hand $\kappa (L) = n.$ Hence \cite{Na85} (see also \cite{Fn07}) applies
and $D = L$ is semi-ample. 
\qed

\section{Almost holomorphic algebraic reductions}
\setcounter{lemma}{0}

We use the same notations as in section 3 and first prove that an almost holomorphic algebraic reduction has in fact a holomorphic model.

\begin{theorem}\label{theorem:holomorphic} Let $X$ be a hyperk\"ahler manifold of dimension $2n$ such that $0 < a(X) < 2n.$
If the algebraic reduction $f$ is almost holomorphic, then $f$ has a Lagrangian holomorphic model.
\end{theorem}

\proof It suffices to show that $L$ is semi-ample.
By Hironaka's flattening theorem \cite{Hi75}
(applied to $\tilde{f} : \tilde{X} \longrightarrow B$) and the normalization (for the resulting source space),
we have an equi-dimensional modification $\hat{f} : \hat{X} \longrightarrow \hat{B}$ of $\tilde{f} : \tilde{X} \longrightarrow B$.
More precisely, there are a normal space $\hat{X}$, a proper bimeromorphic morphism
$$\mu = \hat{\mu} \circ \pi : \hat{X} \longrightarrow \tilde{X} \longrightarrow X,$$
a smooth projective manifold $\hat{B}$,
a birational morphism $\mu_B : \hat{B} \longrightarrow B$ and an equi-dimensional morphism $\hat{f} : \hat{X} \longrightarrow \hat{B}$ 
such that $\mu_B \circ \hat{f} = \tilde{f} \circ \hat{\mu}$ and $\mu_B \circ \hat{f} = f \circ \mu$. Note that we can make $\hat{B}$ smooth as
flatness is preserved under base change. Note also that we can make $\hat{X}$ normal as the normalization map is a finite map. \\
Put $\hat A = \mu_B^*(A)$ so that $\hat A$ is big, nef and semi-ample and set $\hat{L} = \mu^{*}L = \hat{f}^*\hat{A} + \sum a_i\hat{E}_i$, where $\cup \hat{E}_i$
are the exceptional divisors of $\mu$
and $a_i$ are non-negative integers. As $\hat{D}$ and $\hat{f}^*\hat{A}$ 
are Cartier, so is $\sum a_i\hat{E}_i$. Since $f$ is almost holomorphic and $\hat{f}$ is equi-dimensional,
$\hat f(E_i)$ is a divisor on $\hat{B}$.
As $\hat{B}$ is smooth, it is not only a Weil divisor but also Cartier.
Let $C$ be a sufficiently general ample complete intersection curve on $\hat{B}$. Let
$$\hat V = \hat{X} \times_{\hat{B}} C$$ and $\nu : Z \longrightarrow \hat{V} \subset \hat{X}$ be a resolution of
$\hat{V}$. Let $\varphi : Z \longrightarrow C$ be the induced morphism. Let $\eta_Z$ be a K\"ahler class of $Z$. Then, as $\varphi^{*}(\hat{A} \vert C)$ and $\nu^{*}(\sum a_i\hat{E}_i))$ are supported in the fibers of $\varphi$, we have:
$$\int_Z \nu^{*} \hat{L} \wedge \varphi^{*}(\hat{A}\vert C) \wedge \eta_Z^{n-1} = \int_{Z} (\varphi^{*}(\hat{A} \vert C)
+ \nu^{*}(\sum a_i\hat{E}_i)) \wedge \varphi^{*}(\hat{A}\vert C) \wedge \eta_Z^{n-1} = 0\,\, .$$

As $\nu^{*} \hat{L}$  and $\varphi^{*}(\hat{A} \vert C) \in \overline{\mathcal K(Z)}$, they are proportional in $NS(Z)$ by the
Hodge index theorem. So are $\varphi^{*}(\hat{A} \vert C)$ and $\nu^{*}(\sum a_i\hat{E}_i)$.
Consequently
$$\nu^{*}(N(\sum a_i\hat{E}_i)) = \varphi^*(\Theta)$$
for some positive integer $N$ and an effective divisor $\Theta$ on $C$. As $\hat{f}$ is equi-dimensional and $C$ is a general ample complete intersection curve, the Cartier divisor $N(\sum a_i\hat{E}_i)$ on $\hat{X}$ is then of the form $\hat{f}^{*}\Delta$ for some effective Cartier divisor $\Delta$ on $\hat{B}$. Thus, replacing $L$ by some positive
multiple, we have $\mu^{*}\hat L = \hat{f}^{*}(\hat{A} + \Delta)$ for some semi-ample big divisor $\hat{A}$ on $\hat{B}$ and an effective Cartier divisor $\Delta$ on $\hat{B}$. As $\mu^{*}\hat L \in \overline{\mathcal K(\hat{X})}$ (strictly speaking after passing to a
resolution of $\hat{X}$ which does not matter, as we do not need equi-dimensionality any longer), it follows that
$\hat{A} + \Delta \in \overline{\mathcal K(\hat{B})}$. As $\hat{B}$ is projective, this implies that
the divisor $\hat{A} + \Delta$ is nef. As $\hat{A}$ is big and $\Delta$ is effective, the divisor $\hat{A} + \Delta$ is also big. Thus,
$$\kappa(\hat{L}) = \dim\, B = \nu(\hat{A} + \Delta) = \nu(\hat{L})\,\,.$$
As $\hat{L} = \mu^{*}L$, we have $\kappa(L) = \nu(L) > 0$ as well. Thus
$L$ is semi-ample by \cite{Na85} (see also \cite{Fn07}). The morphism given by $\vert mL \vert$ is then Lagrangian fibration by a result of Matsushita \cite{Ma99}. This proves Theorem 4.1. \qed
\vskip .2cm \noindent 
Suppose $0 < a(X) < 2n.$ In order to prove that always $a(X) = n,$ we are reduced to the case that for the general fiber
$F$ has algebraic dimension $a(F) = a(\tilde F) = 0$ and moreover that $\tilde F$ is isotypically semi-simple.
Unfortunately not much is known
about compact K\"ahler manifolds $\tilde F$ with $a(\tilde F) = 0.$ If however $\tilde F$ has a minimal model, things work out:

\begin{proposition} If every isotypically semi-simple compact K\"ahler manifold $Z$  of dimension at most $2n-1$ with $h^0(K_Z) = 1$
has a
minimal model with numerically
trivial canonical bundle (i.e. a bimeromorphically equivalent normal K\"ahler variety which is $\bQ-$Gorenstein with 
numerically trivial canonical class),
then every compact hyperk\"ahler
manifold $X$ of dimension $2n$ has algebraic dimension $a(X) = 0,n,2n.$
\end{proposition}

\proof We must rule out that $1 \leq a(X) \leq n-1.$ We argue by contradiction, hence we are in situation (3.2). We write more precisely
$$ \tilde L = \tilde f^*(A) + \sum_I a_i E_i $$
with $a_i \geq 0.$
Furthermore we have
$$ K_{\tilde X} = \sum_I b_i E_i  $$
with $b_i > 0.$ Let $D_i = E_i \cap \tilde F.$ Then $\tilde L_{\tilde F} = \sum_{I'} a_i D_i$ with $I' \subset I, I' \ne \emptyset$,
 the set of all
$i$ such that $E_i \cap \tilde F \ne \emptyset$ so that by adjunction
$$K_{\tilde F} = \sum_{I'} b_i D_i. \eqno (1)$$
Moreover
$$ \tilde L_{\tilde F} = \sum_{I'} a_iD_i. \eqno (2)  $$
Let $h: \tilde F \rh F'$ be a minimal model of $\tilde F;$ then $K_{F'} \equiv 0.$
Choose a modification $\tau: \hat F \to \tilde F$ from a compact K\"ahler manifold $\hat F$ such that the bimeromorphic map $h: \tilde F
\rh F'$ induces a holomorphic map $\hat h: \hat F \to F'.$ 
Let
$$L' = (\hat h_*(\tau^*(\tilde L_{\tilde F})))^{**}. $$
Since $K_{F'} \equiv 0,$ by (1) every $D_i, i \in I'$ is contracted by $h$. Moreover by (1) and (2)
$$mK_{\tilde F} = \tilde L_{\tilde F} + D'$$
with $D'$ effective and supported in $\bigcup_{I'} D_i.$ Therefore $\tau^*(\tilde L_{\tilde F})$ is on one hand nef, on the other hand
effective with support necessarily in the exceptional locus of $\hat h.$ This is only possible when $\tilde L_{\tilde F} = \sO_{\tilde F}.$
Hence all $a_i = 0$ for $i \in I'$, a contradiction.
\qed
\section{The 4-dimensional Case}
\setcounter{lemma}{0}

In this section we settle Conjecture 1.1 in dimension 4 completely. What still needs to be proved is
\begin{theorem}
Let $X$ be a 4-dimensional hyperk\"ahler manifold. Then $a(X) \ne 1.$
\end{theorem}

\proof Assume to the contrary that $a(X) = 1$. Let $f: X \dasharrow B \simeq \bP_1$ be the algebraic reduction with the setup (3.2);
we set specifically $A = \sO_B(1).$
By (2.4) and (4.1) we know that $a(F) = a(\tilde F) = 0;$ moreover $\tilde F$ is 
isotypically semi-simple. But since $\dim \tilde F = 3$ necessarily $\tilde F$ must be simple. 
We may also assume that $q(\tilde F) = 0;$ otherwise the Albanese map of
$\tilde F$ must be birational onto ${\rm Alb}(\tilde F),$ so that $\tilde F$ has a minimal model, and we conclude by (4.2). 
Since $A = \sO_{\bP}(1),$ we have $h^0(L) = 2;$ we take $F_1, F_2 \in \vert L \vert,$ both necessarily irreducible and set
$$ S = F_1 \cap F_2 $$
as complex spaces. Hence $S$ is a possibly non-reduced complete intersection. Notice that
$$ \pi_*(\sI_E) = \sI_S. \eqno (*) $$
In fact, we have on the level of {\it analytic preimages} (complex subspaces)
$$ \pi^*(S) = \pi^*(F_1) \cap \pi^*(F_2) = (\tilde F_1 + E) \cap (\tilde F_2 + E) = E.$$
In other words $$\pi^*(\sI_S) \cdot \sO_{\tilde X} = \sI_E, $$
where the left hand side denotes the image of $\pi^*(\sI_S) $ in $\sO_{\tilde X}.$ 
Therefore the canonical monomorphism $\sI_S \to \pi_*(\sI_E)$ must be an isomorphism. \\
We first show

\vskip .5cm \noindent
{\bf 5.2 Claim.} {\it $H^q(X,L) = 0$ for $q = 1,3,4$ and $\dim H^2(X,L) = 1.$}

\vskip .5cm \noindent
\proof 
We proceed in several steps.
(1) $H^1(X,L \otimes \sI_S) = 0.$ \\
To verify this vanishing, we deduce from (*)
$$ \pi_*(\tilde f^*(\sO_B(1))) = \pi_*(\sI_E \otimes \tilde L) = \sI_S \otimes L.$$
Thus our claim (1) certainly holds if we can show
$$ H^1(\tilde X,\tilde f^*(\sO_B(1))) = 0.$$
By the Leray spectral sequence (and the projection formula for $\tilde f$), this in turn comes down to
$$ R^1\tilde f_*(\sO_{\tilde X}) = 0. \eqno (**) $$
Since $q(\tilde F) = 0$, the sheaf $ R^1\tilde f_*(\sO_{\tilde X})$ is torsion, supported on a finite set. Thus
if the sheaf would not be $0,$ again the Leray spectral sequence would yield $H^1(\tilde X,\sO_{\tilde X}) \ne 0,$
which is absurd.
\vskip .2cm \noindent
(2) $\chi(L_S) = 0.$ \\
Now there is a constant $K$ such that
$$ a^2 \cdot c_2(X) = K q(a) $$
for all $(1,1)-$classes $a$ (\cite{Fu87}, \cite{GHJ03},23.17). Hence
$$ L^2 \cdot c_2(X) = 0$$
and via Riemann-Roch we obtain
$$ \chi(X,mL) = \chi(\sO_X) = 3 $$
for all $m \in \bZ.$
The Koszul complex
$$ 0 \to L^* \to \sO_X \oplus \sO_X \to L \otimes \sI_S \to 0 $$
gives $\chi(L \otimes \sI_S) = 2 \chi(\sO_X) - \chi(L^*) = 3,$ so that
$$ \chi(L_S) = \chi(L) - \chi(L \otimes \sI_S) = 0.$$
This establishes (2).
\vskip .2cm \noindent
(3) The vanishing (1) and the isomorphism $H^0(\sI_S \otimes L) \to H^0(L) $ gives
$$ H^0(L_S) = 0.$$
\vskip .2cm \noindent Finally we obtain
$$ H^3(L) = H^4(L) = 0 \ ; \ h^2(L) = 1.$$
Concerning $H^2$ we calculate using the adjunction formula $K_S = 2L_S$:
$$ H^2(L_S) = H^0(L_S^* \otimes 2L_S) = H^0(L_S) = 0.$$
Hence by (2):
$$ H^1(L_S) = 0.$$
Therefore $$ H^1(X,L) = 0.$$
Next $$ H^q((X,L) =  0$$
for $q = 4,3$ by Serre duality resp. by a Kodaira vanishing theorem in the K\"ahler case \cite{DP03} (observe $L^2 \ne 0$).
Hence by Riemann-Roch
$$ \dim H^2(X,L) = 1.$$
Thus we completely determined the cohomology of $L$ and Claim (5.2) is established.

\vskip .2cm \noindent
{\it We continue with the proof of Theorem 5.1.} Notice that $h^{1,1}(X) \geq 2,$ otherwise $X$ would be projective. Hence
by \cite{Ca83}, \cite{Fu83-2} there is a smooth hyperk\"ahler deformation $p: \sX \to \Delta$
over the unit disk with the following properties
\begin{itemize}
\item $X \simeq X_0;$
\item there is a line bundle $\sL$ over $\sX$ such that $\sL \vert X_0 \simeq L;$
\item there is a sequence $(t_k)$ in $T$ converging to $0$ such that $X_{t_k} $ is
projective;
\item the set $\Delta_1$ of all $t$ such that $X_t$ is not projective is dense in $\Delta$ with countable complement.  
\end{itemize}
Here $X_t = p^{-1}(t),$ the fiber of $p$ over $t \in \Delta.$ 
Let $L_t = \sL \vert X_t.$
From the knowledge of the cohomology of $L_0,$ semi-continuity and the constancy of $\chi(L_t)$
we obtain immediately (possibly after shrinking $\Delta$) for all $t$:
$$ h^0(L_t) = 2 \ , h^2(L_t) = 1$$
and $$h^q(L_t) = 0 $$
for $q = 1,3,4.$ Therefore
$$ a(X_t) \geq 1 $$
for all $t.$ 
Notice that using (2.7), $a(X_t)$ takes the values $1$, $2$ and $4$; the set $\Delta_2 $ of all $t$ such that $a(X_t) = 1$ is moreover
dense in $\Delta$ with countable complement (otherwise we conclude via a relative algebraic reduction that $a(X_0) \geq 2$).  
We consider the meromorphic map
$$ f_t: X_t \dasharrow B_t \simeq \bP_1$$ defined by $\vert L_t \vert.$ 
Our plan is to apply \cite{AC05} to $X_{t_k}$; \cite{AC05} gives a composition of flops $X_{t_k} \dasharrow X'$ to some other 
projective hyperk\"ahler manifold $X'$
such that
the induced rational map $X' \dasharrow B$ is actually a morphism. But then by \cite{Ma99}, $B_t$ cannot have dimension 1
(a projective hyperk\"ahler 4-fold does not admit a surjective morphism to a curve). 
However in order to be able to apply \cite{AC05} we need to check that
$$ \kappa (F_{t_k}) \leq 0;$$
where $\kappa (F_t)$ is the Kodaira dimension of a desingularisation of a general fiber of $f_t.$ The maps $f_t$
fits together in a family
$$ \mathcal X \dasharrow \bP_1 \times \Delta.$$ 
We introduce a resolution of indeterminacies
$$ \varphi: \tilde {\mathcal X} \to \bP_1 \times \Delta. $$ 
Then we can consider a family $(\tilde F_t)$ of general fibers of $\varphi_t.$ After possibly shrinking $\Delta,$
we may assume that all $\tilde F_t$ are smooth except for $t = 0.$ Here we have an abuse of language and $\tilde F_0$ 
may split in a component which was called $\tilde F_0$ formerly, and possibly some other components. This possible imbellicity
is avoided by considering instead of $X_0$ some $X_s$ with $s  \in \Delta_1, s \ne 0$ and by treating this $X_s$ as our new $X_0.$
Thus we may assume that $p: \tilde X \to \Delta $ is a submersion and that $\tilde F_0$ is smooth. \\  
Now choose a universal number $M$ such that $\vert MK_Z \vert$ defines the Iitaka fibration for all smooth projective
threefolds $Z.$ This number exists by \cite{FM00} and \cite{VZ07} (including references for the general type case and the case
of Kodaira dimension $0$, which actually are not needed here).  
Now by semi-continuity \cite{Gr60} there is a neighborhood $U \subset \Delta$ of $0$ such that 
$$ h^0(MK_{\tilde F_0})  \geq h^0(MK_{\tilde F_{t}}) $$
for all $t \in U.$ Thus 
$$ h^0(MK_{\tilde F_{t_k}}) \leq 1 $$
for all $t_k \in U.$ Therefore $\kappa (\tilde F_{t_k}) = \kappa(F_{t_k}) \leq 0$ for all $t_k \in U.$ 
\qed

\section{Nef line bundles on hyperk\"ahler manifolds}
\setcounter{lemma}{0}

If $X$ is a non-algebraic hyperk\"ahler manifold, then $NS(X)$ is parabolic if and only if $X$ carries a nef non-trivial line bundle $L$, which is 
then unique up to a multiple. We expect that $L$ is actually semi-ample. In this section we give some results pointing in this direction. \\
A line bundle $L$ on a compact complex manifold is {\it hermitian semi-positive} if there exists a (smooth) hermitian metric
on $L$ whose curvature form is semi-positive. Equivalently, there exists a semi-positive $(1,1)-$form $\omega$ such that 
$$ c_1(L) = [\omega]. $$
A hermitian semi-positive line bundle is nef, but the converse is in general not true, see \cite{DPS01}. 

\begin{theorem} Let $X$ be a non-projective hyperk\"ahler manifold of dimension $2n.$ Let $L$ be a non-trivial hermitian semi-positive 
line bundle
on $X.$ Then $a(X) = \kappa(L);$ in particular $\kappa (L) \geq 0.$ 
\end{theorem} 

\proof We use Riemann-Roch in the following form (see e.g. \cite{GHJ03})
$$ \chi(mL) = \sum_{i=0}^n b_i q_X(L)^{i},$$
where $b_i$ are some numbers which do not depend on $L.$
Since $X$ is assumed to be non-algebraic, we have $q_X(L) = 0$ and Riemann-Roch reads
$$ \chi(mL) = b_0 = \chi(\sO_X) = n+1.$$ 
If $h^0(mL) \geq n+1$ for all $m >> 0,$ then $\kappa (L) \geq 1,$ in particular $a(X) \geq 1.$ Since $L$ defines the algebraic
reduction in the sense of (3.2) (recall that $NS(X)$ must be parabolic and that we have only one nef line bundle up to scalars), we obtain
$\kappa (L) = a(X)$. \\
So we may assume that there is a sequence $(m_k)$ converging to $\infty$ and some number $q > 0$ (actually even) such that 
$$ H^q(X,m_kL) \ne 0$$
for all $m_k.$ 
Fix a K\"ahler form $\omega.$
By the Hard Lefschetz Theorem in the semi-positive case \cite{Mo99}, \cite{Ta97}, see also \cite{DPS01}, the
canonical morphism
$$ \wedge \omega^q: H^0(X,\Omega^{2n-q}_X \otimes m_kL) \to H^q(X,m_kL) \eqno (*) $$
is surjective. Thus 
$$ H^0(X,\Omega^{2n-q}_X \otimes m_kL) \ne 0$$
for all $k.$ 
Now we apply \cite{DPS01}, (2.15): one has $a(X) \geq 1$ or $\kappa (L) \geq 0.$ 
In both cases we argue as above and conclude $\kappa(L) = a(X)$. 
\qed

\vskip .2cm \noindent
The arguments of Theorem 6.1 actually sometimes work also in the nef case, namely when the zero locus of a suitable multiplier ideal is not too
large. This leads to the following

\begin{theorem} Let $X$ be a parabolic hyperk\"ahler manifold of dimension $2n \geq 4$. Then $X$ contains a positive dimensional compact subvariety of
dimension at least $2$.
\end{theorem} 

\proof Assume to the contrary that all compact subvarieties of $X$ have dimension at most 1, in particular $a(X) = 0.$ 
Since $NS(X)$ is parabolic, there is, as already mentioned at the beginning of this section, 
a non-trivial nef line bundle $L$, unique up to a scalar . On $L^{\otimes m} $ we introduce a singular metric $h_m$ with multiplier 
ideal $\sI_m$ with  zero locus $V_m.$ We argue similarly as in Theorem 6.1. From Riemann-Roch we deduce the existence of a positive even number $q \geq 2$ 
such that $H^q(X,m_kL) \ne 0$ for a sequence $(m_k)$ converging to $\infty.$ Since $\dim V_{m_k} \leq 1$ by our assumption, we conclude
$$ H^q(X,m_kL \otimes \sI_{m_k}) \ne 0$$
for all $k.$ By the Hard Lefschetz Theorem for nef line bundles \cite{Ta97}, \cite{DPS01}, we obtain the non-vanishing
$$ H^0(X,\Omega^{2n-q}_X \otimes m_kL \otimes \sI_{m_k}) \ne 0.$$ 
Now \cite{DPS01}, (2.15) implies $a(X) \geq 1$ or $\kappa (L) \geq 0.$ Since the only positive-dimensional subvarieties in $X$ are
curves, the first alternative is only possible when $a(X) = 2n-1,$ contradicting (3.6). In the second alternative $X$ contains
a divisor, since $L$ cannot be trivial, again a contradiction. 
\qed

\vskip .2cm
\vskip .2cm
\vskip .4cm \noindent
%\vspace{1cm}
%\small \noindent
Fr\'ed\'eric Campana \\
D\'epartement de Math\'ematiques, Universit\'e de Nancy \\
F-54506 Vandoeuvre-les-Nancy, France \\
frederic.campana@iecn.u-nancy.fr

\vskip .2cm \noindent
Keiji Oguiso \\
Department of Economy, Keio University, Hiyoshi Kouhoku-ku \\
Yokohama, Japan, and Korea Institute for Advanced Study, \\
207-43 Cheonryangni-2dong, Dongdaemun-gu, Seoul 130-722, Korea \\
oguiso@hc.cc.keio.ac.jp

\vskip .2cm \noindent
Thomas Peternell \\
Mathematisches Institut, Universit\" at Bayreuth \\
D-95440 Bayreuth, Germany \\
thomas.peternell@uni-bayreuth.de

\end{document}